\documentclass[12pt,a4paper,fleqn]{article}

\usepackage[hmargin=20mm,top=20mm,bottom=20mm]{geometry}     
                                
\usepackage{pgf,tikz}
\usetikzlibrary{arrows}
\usepackage{caption}
\usepackage{bbm}

\usepackage{mathrsfs}

\usepackage{amsthm}

\usepackage{amsfonts}
\usepackage{amssymb}
\usepackage{amsmath}

\usepackage{centernot}

\usepackage{verbatim}

\usepackage{indentfirst}

\usepackage{amsmath}

\usepackage{color}

\long\def\blue#1{{\color{blue}#1}}

\newtheorem{teo}{Theorem}[section]
\newcommand{\bt}{\begin{teo} }
\newcommand{\et}{\end{teo} }

\newtheorem*{teor}{Theorem}
\newcommand{\bte}{\begin{teor} }
\newcommand{\ete}{\end{teor} }

\newtheorem{prop}[teo]{Proposition}
\newcommand{\bp}{\begin{prop} }
\newcommand{\ep}{\end{prop} }

\newtheorem{lemma}[teo]{Lemma}
\newcommand{\bl}{\begin{lemma} }
\newcommand{\el}{\end{lemma} }

\newtheorem{conj}[teo]{Conjecture}
\newcommand{\bconj}{\begin{conj} }
\newcommand{\econj}{\end{conj} }

\newtheorem{obs}[teo]{Observation}
\newcommand{\bobs}{\begin{obs} }
\newcommand{\eobs}{\end{obs} }

\newtheorem*{lemae}{Lemma}
\newcommand{\ble}{\begin{lemae} }
\newcommand{\ele}{\end{lemae} }

\newtheorem{defi}[teo]{Definition}
\newcommand{\bdefi}{\begin{defi} }
\newcommand{\edefi}{\end{defi} }

\newtheorem{cor}[teo]{Corollary}
\newcommand{\bcor}{\begin{cor} }
\newcommand{\ecor}{\end{cor} }

\newcommand{\bproof}{\begin{proof} }
\newcommand{\eproof}{\end{proof} }

\newcommand{\pic}[1]{  \left\langle #1\right\rangle }

\newcommand{\norm}[1]{  \left\|#1\right\| }
\newcommand{\pare}[1]{  \Bigl(#1\Bigr) }
\newcommand{\paree}[1]{  \left(#1\right) }
\newcommand{\corch}[1]{  \Bigl[#1\Bigr] }
\newcommand{\abs}[1]{  \Bigl|#1\Bigr| }
\newcommand{\llav}[1]{  \Bigl\{#1\Bigr\} }
\newcommand{\llave}[1]{  \left\{#1\right\} }
\newcommand{\abso}[1]{  \left|#1\right| }

 \newcommand{\abse}[1]{  \left|#1\right| }

\newcommand{\st}[2]{\stackrel{#1}{#2}}

\newcommand{\ds}[1]{ \displaystyle{ #1 } }

\newcommand\noi{}
\renewcommand{\varphi}{\gamma}
\renewcommand{\Phi}{F}

\def\nn{\nonumber}
\def\-7{\rightarrow}
\def\RR{\mathbb{R} }

\def\wo{\setminus }

\def\om{\omega}
\def\ZZ{\mathbb{Z} }
\def\EEE{\mathcal{E} }
\def\EE{{ \cal E} }

\def\QQ{{Q} }

\def\PP{{P}}

\def\VV{{\ZZ^d}}

\def\NN{\mathbb{N} }

\def\OOO{\mathcal{O} }
\def\XXX{\mathcal{X} }
\def\11{\mbox{\normalfont\bfseries{1}} }

\def\FFF{\mathcal{F} }
\def\GGG{\mathcal{G} }
\def\LLL{{L} }

\def\WWW{{W} }

\def\cleq{\preceq}
\def\cgeq{\succeq}

\def\Ber{{B} }
\def\rbeta{{\rho}}
\def\aa{{a}}
\def\bb{{b}}

\def\Hem{\mbox{Hem}}

\def\S{{S}}
\def\V{{V}}
\def\insv{{}}

\title{Phase transition for the dilute clock model}
\author{In\'es Armend\'ariz\thanks{Departamento de Matem\'atica, UBA, Buenos Aires, Argentina. Email:\,{\tt iarmend@dm.uba.ar}},
Pablo A. Ferrari\thanks{Departamento de Matem\'atica, UBA and IMAS-CONICET, Buenos Aires, Argentina, and IME-USP, S\~ao Paulo, Brazil. Email:\,{\tt pferrari@dm.uba.ar}},
Nahuel Soprano-Loto\thanks{IMAS-CONICET, Buenos Aires, Argentina. Email:\,{\tt nsloto@dm.uba.ar}}}
\begin{document}
\sloppy
\maketitle
 
\begin{abstract}
\noindent 
We prove that phase transition occurs in the dilute ferromagnetic nearest-neighbour $q$-state clock model in $\ZZ^d$, for every $q\geq 2$ and $d\geq 2$.
This follows from the fact that the Edwards-Sokal random-cluster representation of the clock model stochastically dominates a supercritical Bernoulli bond percolation probability,  a technique that has been applied to show phase transition for the low-temperature Potts model.
The domination involves a combinatorial lemma which is one of the main points of this article.
\\
\\
{\em AMS 2000 Mathematics Subject Classification}: 82B\\
\\
{\em Keywords}: Dilute clock model, phase transition.

\end{abstract}


\section{Introduction}
\parskip 2mm
The $q$-state clock model assigns a random spin to each site of $\ZZ^d$. The spins take values in a discrete set $S$ of equidistant
angles or hours, hence the name.
Let $\sigma=(\sigma_x, x\in\ZZ^d)$ be a spin
configuration, $\sigma_x$  the angle of the spin at $x\in\ZZ^d$.
Let $\EEE(\ZZ^d):= \{\langle xy \rangle: \|x-y\|=1\}$ be the set of edges connecting nearest neighbour sites, $\norm{\cdot}$ the Euclidean norm.
We study the dilute clock model associated to a disorder, namely a collection 
\begin{equation}
\label{disorder}
J=(J_{\langle xy \rangle}: \langle{xy}\rangle \in \EEE(\ZZ^d))
\end{equation} 
of independent identically distributed Bernoulli random variables with parameter $p$. 
A disorder realization $J$ and a finite set $\Lambda\subset \ZZ^d$ determine the  Hamiltonian on spin configurations:
\begin{align}
\label{ham0}
H_{\Lambda,J}(\sigma):=\sum_{\substack{\langle{xy}\rangle\in \EEE(\ZZ^d)\\\llave{x,y}\cap \Lambda\neq \emptyset }}J_{\langle xy\rangle}\big(1-\cos(\sigma_x-\sigma_y)\big).
\end{align}
When $q=2$, we recover the Ising model; as $q\to \infty$, the clock model approximates the $XY$ model, which has a continuum of spin angles.

Given a  set $\Lambda\subset \ZZ^d$ and configurations $\sigma,\eta\in S^{\ZZ^d}$, we write
\begin{equation}
\label{bc}
\sigma\st{\Lambda}{=}\eta\quad \mbox{if}\quad \sigma_x=\eta_x\ \forall \,x\in \Lambda.
\end{equation} 
The configuration $\eta$ plays the role of a boundary condition.
The specification $\mu_{\Lambda,J}^{\eta}$ associated to a finite set
$\Lambda$, a disorder $J$, and
 a boundary condition $\eta$ is the probability 
\begin{align}
\label{specif}
  \mu^\eta_{\Lambda,J} ( \sigma ) :=\frac{1}{Z^\eta_{\Lambda,J}}e^{-\beta
    H_{\Lambda, J}(\sigma)} \11 [\sigma\st{\Lambda^c}{=}\eta]\,,
\end{align}
where $\beta>0$ is a parameter proportional to the inverse temperature, $Z^\eta_{\Lambda,J}$ is the normalizing constant and $\Lambda^c=\ZZ^d\setminus\Lambda$.
A Gibbs measure associated to the disorder $J$ is a probability $\mu_J$ that satisfies the DLR condition:
\begin{align}
\mu_Jf=\int_{S^{\ZZ^d}}\mu_J(d\eta)\, \mu_{\Lambda,J}^\eta f\,
\end{align}
for every finite subset $\Lambda\subset \VV$ and every local
function $f:S^{\ZZ^d}\-7 \RR$.
Here, $\mu f$ denotes the
expectation of $f$ with respect to $\mu$.
The underlying $\sigma$-algebra where the Gibbs measures and the specifications are defined is the one generated by projections over finite subsets of $\ZZ^d$.
We call $\GGG_J$ the set of Gibbs measures associates to $J$.
Since $S$ is finite,  $\GGG_J$ is not empty.
In case $\abse{\GGG_J}>1$, we say that phase co-existence occurs.

The homogeneous version of the model is obtained by taking $p=1$ or, equivalently, $J_{\pic{xy}}\equiv 1$ for every $\pic{xy}$.
In this case, 
non-uniqueness methods such as the Pirogov-Sinai theory
\cite{PS75} or reflection positivity 
as in Fr\"olich, Israel, Lieb and Simon 
\cite{FILS78}, see also Biskup \cite{biskup},
prove that, for sufficiently low temperature, there exist at least $q$ different Gibbs measures.
On the other hand, when the temperature is large enough, techniques similar to those developed by Dobrushin
\cite{dobrushin} or by van den Berg and Maes \cite{vdbm} show that there exists only one Gibbs measure. Phase transition occurs when a system undergoes
a change in its phase diagram depending on the value of a parameter;
these results hence establish occurrence of phase transition for the homogeneous clock model.

Both Pirogov-Sinai theory and reflection positivity depend on the graph determined by the interacting edges in the Hamiltonian \eqref{ham0} being symmetric, an assumption that breaks down for the properly dilute model $p<1$.
Instead, our main tool is the Fortuin-Kasteleyn
random-cluster representation \cite{KF69}, originally introduced for the Ising and the Ashkin-Teller-Potts models, and then generalized to arbitrary models by Edwards and Sokal \cite{es}; here, we build the clock model random-cluster representation in detail. The core idea of this approach is to relate non-uniqueness of Gibbs measure in the  
statistical-mechanical model to the existence of an 
infinite cluster in the random-cluster model: a percolation problem. It was first applied by Aizenman, Chayes, Chayes and Newman  to study 
the phase diagram of the dilute Ising and Potts models in \cite{ACCN87}; we presently adapt their ideas to our context.

Precisely, we derive a lower bound for the critical temperature: for every dimension $d$ and every
number $q$ of spins, we take $p$ sufficiently large to guarantee that the disorder 
almost surely contains an infinite bond-percolation cluster, and then determine a value $\beta_0=\beta_0(q,d,p)>0$ such that there is more
than one Gibbs measure at inverse temperatures $\beta >\beta_0$, for almost all disorders $J$.
A crucial step in the proof consists of
dominating the random-cluster probability associated to the clock model by
a supercritical Bernoulli product probability on the bonds. While
for the Potts model this domination is immediate, the clock model requires a
combinatorial argument, given in Lemma \ref{combinatorics}. 

The following is our main result.

\begin{teo}\label{diluted}
Let $p>p_c$, where $p_c$ is the Bernoulli bond percolation critical probability
in $\ZZ^d$.
Let the disorder $J$ be distributed as a product $\PP_p$ of i.i.d. Bernoulli random variables with parameter $p$. 
Then
there exists $\beta_0>0$ such that, for $\beta>\beta_0$, the $q$-state dilute clock model associated to the random specifications $\mu^\eta_{\Lambda,J}$
defined in \eqref{specif} exhibits phase
co-existence for $\PP_p$-almost every realization of $J$. More precisely,
$\beta>\beta_0$ implies $\PP_p(J: \abse{\GGG_J}\geq q)=1$.
\end{teo} 

The value $\beta_0$ in the later theorem depends on $d$,
$q$ and $p$. 
We show in the Appendix that, for fixed $d$ and $p$, $\beta_0(q,d,p)\sim q^2\log
q$ as $q\-7\infty$, the same asymptotics provided by Pigorov-Sinai theory and reflection positivity in the $2$-dimensional  homogeneous case.
In particular, 
$\lim_{q\to\infty} \beta_0(q,d,p) =\infty$, implying that our approach
is not suitable to study the $XY$ model;
see van Enter, K\"ulske and Opoku \cite{vEKO11} for results concerning the approximation of the $XY$ model via the clock model.
On the other hand, for $d\geq 3$ and $p=1$, reflection positivity computes a threshold $\beta_0$ 
independent of $q$, 
see Maes and Shlosman \cite{MS11} for a discussion.

The ideas presented in this article can be further developed in two directions,
which are explored by Soprano-Loto
in
collaboration with Roberto Fern\'andez in a separate article \cite{F-SL}.
The first one is
a generalization of the current work to the so called Abelian spin models; see Dub\'edat \cite{D11} for a precise definition.
The second direction of research seeks to obtain a uniqueness criterion, also via  random-cluster representation, at a higher level of generality. 





\noindent {\bf Organization of the article.} We introduce the random-cluster model and state the results leading to the proof of Theorem \ref{diluted} in Section \ref{aux}. Section \ref{proofs} contains the proofs and the Appendix collects some auxiliary computations.

\section{Clock model and random-cluster in a finite graph}
\label{aux}

We define the clock model and its random-cluster representation for a fixed  non-oriented finite graph $(\V,\EEE)$
without loops or multiple edges, and not necessarily connected.
We fix a non-empty subset $U\subset V$ playing the role of \emph{boundary}.
For simplicity, we suppose there are no edges connecting vertices in $U$:
$\llave{\pic{xy}\in \EEE:\llave{x,y}\subset U}=\emptyset$. 

In the case of the dilute clock in a finite set $\Lambda\subset \ZZ^d$, the boundary is given by 
$\partial\Lambda:= \{y\in\ZZ^d\setminus\Lambda,\, \exists x\in \Lambda: \|x-y\|^2= 1\}$, and
the vertex and edge sets are
 \begin{equation}
\label{dilmodel}
\Lambda\cup \partial\Lambda \quad\mbox{ and } \quad \{\pic{xy},\,\{x,y\} \not\subset \Lambda^c,\, \|x-y\|=1,\,J_{\langle xy\rangle}=1\}.
 \end{equation}

\vspace{-0.7cm}
\paragraph{The clock model.}
Let $\S$ be the set of angles defined by
\begin{equation}
\label{spinset}
\S:=\llave{\frac{2\pi i}q :i=0,\ldots,q-1}. 
\end{equation}
Elements of $\S$ are called \emph{spins} and denoted  $a$,  $b$ and $c$, while \emph{spin} or \emph{vertex-configurations} in $\S^\V$ are  denoted  by $\sigma$ and $\eta$. 

The clock Hamiltonian $H=H(V,\EEE)$ is the function $H:S^V\-7\RR$  defined by
\begin{align}
H(\sigma):= \sum_{\pic{xy}\in\EEE}\big(1-\cos(\sigma_x-\sigma_y)\big).
\end{align}
We write $\sigma\st{U}{=}a$ when $\sigma_x=a$ for all $x\in U$. 
The clock probability $\mu=\mu(V,U,\EEE, \beta)$ with $0$-boundary condition is defined as 
\begin{align}\label{see}
\mu(\sigma) := \frac{1}{Z}e^{-\beta H(\sigma)}\11 [\sigma\st{U}{=}0],
\end{align}
where $\beta$ is a strictly positive parameter and $Z=Z(V,U,\EEE,\beta)$ is the normalizing constant. 
\vspace{-0.4cm}
\paragraph{The random-cluster measure.}
Define a weight function $ \WWW: S\-7(0,1]$  by
\begin{equation}
\label{wst}
   \WWW(\aa):= e^{-\beta (1-\cos a)}
\end{equation}
and let ${\cal I}:=\{\WWW(\aa),\,\aa \in \S\}$ be its image.
This set has cardinality $|{\cal I}|=k+1$, where $k= q/2$ for even $q$ and  $k=(q-1)/2$ for odd $q$.  Write ${\cal I}=\{t_0,t_1,\dots, t_k\}$ with $0<t_0<t_1<\dots<t_k=W(0)=e^{-\beta (1-\cos0)}=1$, and denote  
\begin{equation}
\label{ewights}
r_0:= t_0, \quad 
r_i:= t_i- t_{i-1}, \quad   1\le i\le k.
\end{equation}
By construction $0\le r_i\le 1$ for all $0\le i\le k$, and $\sum_i r_i=1$. 

Let $\theta$ be the probability on ${\cal I}$ given by
\begin{equation}
\label{hatfi1}
\theta(t_i):= r_i, \quad 0\le i\le k,
\end{equation}
and let  $\hat\phi=\hat\phi(\EEE,\beta)$
be the product measure on the set of edge-configurations $\omega\in{\cal I}^\EEE$ with marginals~$\theta$:
\begin{equation}
\label{hatfi2}
\hat\phi(\omega)
:=\prod_{\pic{xy}\in\EEE}\theta(\omega_{\pic{xy}}).
\end{equation}
We say that an edge-configuration $\om\in {\cal I }^{\EEE}$ and a vertex-configuration
$\sigma\in S^V$ are \emph{compatible}, and write $\om \cleq \sigma$, if the value of $\om$ on any edge is dominated by the weight of the gradient of $\sigma$ over that edge: 
\begin{align}
\label{compat}
\om\cleq\sigma\quad\Leftrightarrow \quad \om_{\pic{xy}}\leq \WWW(\sigma_x-\sigma_y) \mbox{ for every } \pic{xy}\in \EE.
\end{align}
Notice that  if $\om \cleq \sigma$ and $\om_{\pic{xy}}=1$, then
$\sigma_x=\sigma_y$; on the other hand, 
$\om_{\pic{xy}}=0$ imposes no restriction on the values of $\sigma_x$ and $\sigma_y$.
 
We define the \emph{random-cluster probability} $\phi=\phi(V,U,\EEE,\beta)$ on ${\cal I}^{\EEE}$ as the
measure obtained from
$\hat\phi$ by assigning to each edge-configuration $\om$ a weight proportional to the number of vertex-configurations 
$\sigma$ that are compatible with $\om$ and satisfy the boundary condition, using $\hat\phi$ as reference measure:
\begin{equation}
\label{perro}
 \phi(\om)
:=\frac{1}{Z } \;
    \big|\{\sigma\insv :\sigma\cgeq
        \om,\,\sigma\st{U}{=}0\}\big|\; \hat\phi(\om).
\end{equation}
 Here $Z$ is the same normalizing constant appearing in \eqref{see}.


\vspace{-0.4cm}
\paragraph{The Edwards-Sokal coupling.} 

Let $\hat\mu=\hat\mu(\V,U)$ be the uniform probability on the set  of vertex  configurations $\S^V$
that are identically 0 at sites in $U$: 
\begin{align}
\hat\mu(\sigma) := \frac{1}{q^{\abse{V\wo U}}}\11 [ \sigma\st{U}{=}0 ].
\end{align}
We define a joint edge-vertex probability $\QQ=\QQ(V,U,\EEE, \beta)$ on ${\cal I}^\EEE\times S^V$  by
\begin{align}
\label{Q}
  \QQ(\om,\sigma):=
  \frac{1}{ Z'}\  \11 [\om\cleq \sigma]\ \hat\phi(\om)\  \hat\mu(\sigma),
\end{align}
where  $Z':= Z/q^{|V\setminus U|}$ with $Z$ as in \eqref{see}.
That is, $\QQ$  is the product
probability $\hat\phi\times\hat\mu$ conditioned to the compatibility event
$\llave{(\om,\sigma): \om\cleq\sigma}\subset {\cal I}^\EEE\times S^V$. 

\begin{teo}\label{edso}
Edwards-Sokal \cite{es}.

\noindent The measures $\phi$ and $\mu$ are respectively the first and second marginals of $Q$.
\end{teo}

We prove this theorem in 
Section \ref{proofs}. The measure $\QQ$ can be seen as a coupling between the clock measure $\mu$ and the random-cluster measure $\phi$. 
As a corollary, it follows that  the conditional
distribution under $Q$ of $\sigma$ given $\omega$ is uniform
on the set of configurations compatible with $\omega$  and such that $\sigma\st{U}{=}0$:
\begin{align}
  \label{csssa}
Q\big(\sigma\,|\,\omega\big) =
\frac{\QQ(\om,\sigma)}{\sum_{\sigma'}\QQ (\om,\sigma') } 
\;=\; 
\frac{\hat\mu(\sigma)\;\11 [\om\cleq\sigma]}{\hat\mu\big(\sigma' \insv : \om\cleq\sigma'\big)}\,.
\end{align}
This implies that a random vertex-configuration distributed according to  $\mu$ may be sampled by first choosing
an edge-configuration $\om$ with law $\phi$, and then sampling a vertex-configuration uniformly among those that are compatible with $\om$ and satisfy the boundary restriction.
That is, 
\begin{align}\label{sssa}
\mu(\sigma)=\sum_{\om\in {\cal I}^\EEE}\frac{\hat\mu(\sigma)\,\11 [\om\cleq\sigma]}{\hat\mu\big(\sigma' : \om\cleq\sigma'\big)}\; \phi(\om).
\end{align}

Given $x,y\in V$ and $\om\in {\cal I}^\EEE$, we denote
$x\st{\om}{\longleftrightarrow}y$ if there is a path of vertices $x_1,\dots, x_n\in V$ with $x_1=x$,
$x_n=y$,  $\pic{x_ix_{i+1}}\in \EEE$ and $\om_{\pic{x_ix_{i+1}}} =1$ for $1\le i\le n-1$. We say that $x$ is \emph{connected} to $U$ by an $\omega$-open path,
and write $x\st{\om}{\longleftrightarrow}U$, when $x\st{\om}{\longleftrightarrow}y$ for some $y\in U$;
let $x\st{\om}{\centernot\longleftrightarrow} U$ denote the complementary event.
The $\mu$-marginal of the spin at $x$ can be related to the connection probabilities between $x$ and the boundary, under $\phi$ and $Q$: 
\begin{align}
  \label{i14}
\mu(\sigma:\sigma_x=a) \;=\; \phi\big(\omega: x\st{\om}{\longleftrightarrow} U\big) \11[a=0] \;+\; Q\big((\om, \sigma): \sigma_x=a,\, x\st{\om}{\centernot\longleftrightarrow} U\big).
\end{align}
Identity \eqref{i14} follows immediately from the coupling of Theorem \ref{edso} and the inclusion $\{(\omega,\sigma): x\st{\om}{\longleftrightarrow} U\}\subset\{(\omega,\sigma):\sigma_x=0\}$.

The coupling of Theorem \ref{edso} also implies that the $\mu$-probability of seeing a 0 at any site $x$ is larger than the probability of seeing any other spin plus the $\phi$-probability that $x$ be connected to the boundary.
This is the content of the next result; its proof depends crucially on the combinatorial Lemma \ref{combinatorics} stated later.  

\begin{prop}
\label{correlations} Positive correlations.

\noindent For any vertex $x\in \V$ and any spin $a\neq 0$,
\begin{align}
  \label{i15}
\mu(\sigma:\sigma_x=0)\;\ge\; \mu(\sigma:\sigma_x=a)\,+ \,\phi\big(\omega: x\st{\om}{\longleftrightarrow} U\big).
\end{align}
\end{prop}

\paragraph{Stochastic domination.}
Given $ I\subset \RR$, consider the partial order on $I^\EE$ defined by  $\om \leq \om'$ if and only if $\om_{\pic{xy}}\leq
\om'_{\pic{xy}}$ for every $\pic{xy}\in \EE$.
A function $f: I^{\EE}\-7 \RR$ is said to be increasing if $f(\om)\leq f(\om')$ whenever $\om\leq \om'$, while
an event $E\subset I^{\EE}$ is said to be increasing when its indicator function $f(\om) = \11[\om\in E]$ is. 
Given two probabilities $P$ and $P'$ on $ I^{\EE}$, we say that $P$ is stochastically dominated by $P'$,
and write $P\leq_{st}P'$, if and only if $Pf\leq P'f$ for every increasing
$f: I^{\EE}\-7 \RR$.
This is equivalent to $P(E)\le P'(E)$ for any increasing event $E$.

Given $ \rho \in [0,1]$, let $B_\rho$ be the Bernoulli product measure on   $\{  0,1  \}^\EEE$  with parameter $\rho$
.
In order to  stochastically compare $\phi$ and $B_\rho$
we consider them defined on   
the common space $I^\EEE$, where $I=\{0\}\cup {\cal I}$.


\begin{teo}\label{domination}
Stochastic domination.

\noindent For any $\rho\in [0,1)$ there exists $\beta_0=\beta_0(\rho)>0$, independent of the
graph $(V,\EEE)$ and the boundary $U$, such that, if $\beta\geq\beta_0$, $B_\rho$ is
stochastically dominated by $\phi$.
\end{teo}

The key to the proofs of Proposition \ref{correlations} and Theorem \ref{domination} is the
following combinatorial lemma, proved  in Section \ref{proofs}.

\begin{lemma}\label{combinatorics}
For every $x\in V$, $a\in S$ and $\om\in {\cal I}^\EEE$,
\begin{align}\label{ui9}
\big| \{ \sigma\insv :\sigma\cleq \om, \sigma\st{U}{=}0,\sigma_x=a \}\big|
\;\leq\;
\big| \{ \sigma\insv :\sigma\cleq \om, \sigma\st{U}{=}0,\sigma_x=0 \}\big|.
\end{align}
Equivalently,
\begin{align}
\label{ui99}
  \hat\mu(\sigma: \sigma_x=a,\,\sigma\cleq \om) \le  \hat\mu(\sigma: \sigma_x=0,\,\sigma\cleq \om).
\end{align}
\end{lemma}

The lemma in fact holds for any spin set $\S'$ and weight function $W'$ provided they satisfy certain symmetry properties: for any pair of elements $a, b \in \S'$ it must be possible to define a reflection $R=R_{a,b}:\S'\to\S'$, $R(a)=b$, such that  {\it i)} it splits $S'$ into two hemispheres $\mbox{Hem}(a)$ and Hem$(b)$, $a \in \mbox{Hem}(a)$, $b \in \mbox{Hem}(b)$, in such a way that $W'(c-R(d))<W'(c-d)$ implies $c$ and $d$ belong to the same hemisphere, and {\it ii)} $R$ preserves the compatibility of neighbouring vertices when applied to both spins. These extensions are explored in detail in \cite{F-SL}.

In the dilute Potts model with $q$ spins, the
Hamiltonian is given by $\sum_{\pic{xy}}J_{\pic{xy}} \11 [\sigma_x\neq\sigma_y]$,
and the associated random-cluster probability is defined on $\llave{0,1}^{\EE}$;
see \cite{GHM01,G06}, for example.
In this case, if  $\sigma$ and $\om$ are compatible,  the values of $\sigma_x$ and $\sigma_y$ must coincide whenever $\om_{\pic{xy}}=1$, and there are no restrictions if $\om_{\pic{xy}}=0$.
Call a connected component  of the graph $(\V, \{\pic{xy}:\om_{\pic{xy}}=1\})$ an $\om$-cluster.
Then $\om\cleq \sigma$
implies that $\sigma$ is constant over each of the $\om$-clusters and the values  achieved on different clusters  not connected with $U$ can take any value in $\{1,\dots,q\}$. 
Hence,  for the diluted Potts model, the combinatorial term appearing in expression \eqref{perro} reduces to
\begin{align}
| \{\sigma\insv :\sigma\cgeq
        \om,\sigma\st{U}{=}0\} |= q^{\textnormal{number of $\om$-clusters}}.
\end{align}
In contrast, for the clock model, the larger range of  edge-configurations in ${\cal I}^\EE$ gives rise to a more
delicate combinatorial structure which will be managed using the inequality \eqref{ui9}.


\section{Proofs}\label{proofs}

\paragraph{Proof of Theorem \ref{diluted}: Phase co-existence.}
Let us identify a disorder $J$ defined in \eqref{disorder} with its associated set of open edges 
\begin{align}\label{t}
\llave{\langle x y\rangle \in \EEE(\ZZ^d): J_{\langle xy\rangle }=1}.
\end{align}
We say that $C\subset \ZZ^d$ is a $J$-open cluster if it is a maximal set with the property that $x\st{J}{\longleftrightarrow}y$ for all $x,\,y \in C$.
Denote $x\st{J}{\longleftrightarrow}\infty$ when $x$ belongs to an infinite $J$-open cluster. Let $p_c$ be the critical value for bond percolation in $\ZZ^d$. If $p>p_c$ then $P_p(J:x\st{J}{\longleftrightarrow}\infty)>0$; see \cite{GHM01,G06} and references therein  for a treatment of percolation theory.   

Let $\rho \in (0,1)$ be such that $p\rho>p_c$ and let $J'$ be an independently sampled $P_\rho$-disorder. 
Denote by $JJ'$ the set of vertices that are open for both $J$ and $J'$, note that $JJ'$ is a $P_{p\rho}$-disorder. Also, once $J$ is fixed, $JJ'$ 
is a random thinning, each open edge of $J$ is kept open with probability $\rho$ and closed with probability $(1-\rho)$, independently. 

Let $\XXX\subset \llave{0,1}^{\EEE(\ZZ^d)}$ be the set of disorders $J$ such that there is an infinite $JJ'$-open cluster with probability $1$:
\begin{align}
  \label{dxxx}
\XXX:= \{J:P_\rho(J': \mbox{ there is an infinite $JJ'$-open cluster})=1\}.
\end{align}
From Fubini's Theorem, the fact that $JJ'$ is a $P_{p\rho }$-disorder, and $p\rho>p_c$, it is easy to see that $P_p(\XXX)=1$. 
Also,
\[\big\{J':\mbox{ there is an infinite $JJ'$-open cluster}\big\}=\bigcup_{x\in\ZZ^d}\big\{J': x\st{JJ'}{\longleftrightarrow}\infty\big\}.
\]
Hence, for each $J\in \XXX$, there exists a vertex  $x\in\ZZ^d$ belonging to an infinite $JJ'$-open cluster with positive $P_\rho$-probability:
\begin{align}
\label{l}
 P_{\rho}\big(J':x\st{JJ'}{\longleftrightarrow}\infty\big)>0.
\end{align}

Let $\beta_0=\beta_0(\rho)$ be as in the statement of Theorem \ref{domination}. 
Fix  a disorder $J\in \XXX$ and a vertex $x$ satisfying \eqref{l}. 
Given $n\in \NN$, let $\Lambda_n:=[-n,n]^d\cap \ZZ^d$ and consider the choices 
\begin{align}
\label{choice}
\V=\Lambda_n\cup \partial \Lambda_n,\quad \EEE=\{\pic{xy}\in\EEE(\ZZ^d):\{x,y\}\cap \Lambda_n\neq \emptyset,\; J_{\langle xy\rangle}=1\},\quad  U=\partial \Lambda_n,
\end{align}
for the vertex, edge and boundary sets in Section \ref{aux}. 
Let $\mu$, $\phi$ and $B_\rho$ denote the  clock probability on $S^{V}$, random-cluster distribution on ${\cal I}^\EEE$ and product Bernoulli probability on 
$\{0,1\}^{\EEE}$ associated to this choice, respectively. Note that $\mu=\mu_{\Lambda_n, J}^{0}$ as defined in \eqref{specif} with the convention that  the superscript $a$  in $\mu_{\Lambda_n, J}^a$ indicates the boundary condition 
$\eta_y\equiv a$ on $\partial \Lambda_n$.

Since the event $\{x\st{\om}{\longleftrightarrow}U\}$ is increasing, Theorem \ref{domination} implies
\begin{align}
\phi\big(\om: x\st{\om}{\longleftrightarrow}U\big)
\geq
B_\rho \big(\om: x\st{\om}{\longleftrightarrow}U\big)=
P_{\rho}\big(J' : x\st{JJ'}{\longleftrightarrow}\partial\Lambda_n\big)
\geq 
P_{\rho}\big(J' : x\st{JJ'}{\longleftrightarrow}\infty\big). \notag
\end{align}
Replacing in  \eqref{i15} with $\mu=\mu_{\Lambda_n,J}^0$, we obtain
\begin{align}
\mu_{\Lambda_n,J}^0(\sigma:\sigma_x=0)\,\ge\, \mu_{\Lambda_n, J}^0(\sigma:\sigma_x=a)+P_{\rho}\big(J' : x\st{JJ'}{\longleftrightarrow}\infty\big),\qquad \hbox{ for any }a\neq 0.
\end{align}
We conclude that any weak limit $\mu_J^0$ of $\mu_{\Lambda_n, J}^0$ as $n\to \infty$ will satisfy 
\begin{align}
\label{limit}
\mu_J^0(\sigma:\sigma_x=0)>\mu_J^0(\sigma: \sigma_x=a)\qquad \hbox{ for any }a\neq 0.
\end{align}
By the rotational symmetry in the set $S$ of spins, the same holds with any boundary condition $b$: the weak limit $\mu_J^b$ assigns maximal probability to having a spin $b$ at $x$, 
$\mu_J^b(\sigma: \sigma_x=b)>\mu_J^b(\sigma: \sigma_x=a), \ a\neq b$, and therefore the $q$-Gibbs measures $\mu_J^b,\,b\in S$, must be different.
\qed


\paragraph{Proof of Proposition \ref{correlations}: Positive correlations.}


For any spin $a\neq 0$, by \eqref{sssa} and the fact that $x\st{\om} {\longleftrightarrow}U$ implies $\sigma(x)=0$, 
\begin{align}
&\mu(\sigma:\sigma_x=a)
\;=\;
\sum_{\om:\,x\st{\om}{\centernot\longleftrightarrow} U} \frac{\hat\mu(\sigma:\sigma_x=a,\om\cleq\sigma)}{\hat\mu(\sigma:\om\cleq\sigma)}\;\phi(\om)  \notag\\[3mm]
&\qquad\le\; 
\sum_{\om:\,x\st{\om}{\centernot\longleftrightarrow} U} \frac{\hat\mu(\sigma:\sigma_x=0,\om\cleq\sigma)}{\hat\mu(\sigma:\om\cleq\sigma)}\;\phi(\om)
\;=\;Q\big((\om, \sigma): \sigma_x=0,\, x\st{\om}{\centernot\longleftrightarrow} U\big), \notag
\end{align}
where the inequality holds by \eqref{ui99}. Apply \eqref{i14} to conclude. \qed

\paragraph{Proof of Theorem \ref{domination}: Stochastic domination.}
The measure $\phi$ gives positive probability to every edge configuration. Under this hypothesis,  Holley's inequality (Theorem 4.8 of \cite{GHM01} for instance),
asserts that the stochastic domination
$\Ber_{\rho} \le_{st} \phi$ follows from the single-bond inequalities
\begin{align}
\label{199}
\rho \le \phi\big(\om: \om_{\pic{xy}}=1\big|\,\om:\om\st{\EE\setminus\langle xy \rangle}{=} \om'\big)=:\alpha(\langle xy\rangle,\om'), \qquad \pic{xy}\in \EE,\;\;\om'\in  {\cal I}^{ \EE}.
\end{align}
Given $t\in {\cal I}$, we define $t_{\pic{xy}}\om'\in {\cal I}^{\EE}$ by
\begin{align}
  \label{txy}
(t_{\pic{xy}}\om')_{\pic{xy}}=t\qquad\hbox{and}\qquad
t_{\pic{xy}}\om' \st{\EE\setminus\langle xy \rangle}{=} \om'\notag .
\end{align}
Omitting the dependence of $\alpha$ on $(\langle xy\rangle,\om')$ in the notation, 
\begin{align}
\alpha=  \frac{\phi(1_{\pic{xy}}\om')}{\sum_{i=0}^k\phi((t_i)_{\pic{xy}}\om')}
  =\frac{r_k |\{\sigma\insv : \sigma\cgeq 1_{\pic{xy}}
        \om',\sigma\st{U}{=}0\}|}
  {\sum_{i=0}^{k}r_i|\{\sigma\insv : \sigma\cgeq (t_i)_{\pic{xy}}
        \om',\sigma\st{U}{=}0\}|},
\end{align}
and
\begin{align}
\alpha^{-1}= 
\sum_{i=0}^{k}\frac{r_i}{r_k} \ \frac{|\{\sigma\insv :
      \sigma\cgeq (t_i)_{\pic{xy}}
      \om',\sigma\st{U}{=}0\}|}{|\{\sigma\insv :
      \sigma\cgeq 1_{\pic{xy}}
      \om',\sigma\st{U}{=}0\}|}.\label{mapr} 
\end{align}
Let $(V,\tilde{\EE})$ be the auxiliary graph obtained from $(V,\EEE)$ by adding all edges connecting vertices in $U$: 
\begin{align}
\tilde\EE:= \EE\cup \big\{\pic{uv}:\,\{u,v\}\subset U\big\}.
\end{align}
Let 
$\tilde\om\in {\cal I}^{\tilde{\EE}}$ be defined by
\begin{align}
\tilde\om\st{\tilde\EE\setminus\EE}{=}1\qquad \hbox{ and }\qquad
\tilde\om\st{\EE}{=}\om'\notag .
\end{align}
Extend the definition of $t_{\pic{xy}}\tilde\om \in {\cal I}^{\tilde{\EEE}}$ and the compatibility notion $\sigma \cleq \tilde\om$ to the enlarged graph in the obvious way and use the rotation invariance of $S$ to get 
\begin{align}
|\{\sigma\insv : \sigma\cgeq (t_i)_{\pic{xy}} \om',\sigma\st{U}{=}0\}|=\frac{1}{q}|\{\sigma\insv: \sigma\cgeq (t_i)_{\pic{xy}} \tilde\om\}|,
\end{align}
and replacing in (\ref{mapr}),
\begin{align}
\label{eert}
\alpha^{-1}=\sum_{i=0}^{k}\frac{r_i}{r_k}\frac{|\{\sigma\insv: \sigma\cgeq (t_i)_{\pic{xy}} \tilde\om\}|}{|\{\sigma\insv: \sigma\cgeq 1_{\pic{xy}} \tilde\om\}|}.
\end{align}
For $0\leq i\leq k$, let
\begin{align}
\label{ksubi}
K_i:= |\{(\aa,\bb)\in  S\times S: \WWW(\aa-\bb)= t_i\}|.
\end{align} 
We have
\begin{align}
  |\{\sigma\insv: \sigma\cgeq (t_i)_{\pic{xy}} \tilde\om\}|
 =&\sum_{j=i}^{k}|\{\sigma\insv : \sigma\cgeq
      \tilde\om, \WWW(\sigma_y-\sigma_x)= t_j\}| \\[0.5cm]
  =&\sum_{j=i}^{k}K_j\, |\{\sigma\insv : \sigma\cgeq \tilde\om,
      \sigma_y=0,\sigma_x=\aa_j\}|,
\end{align}
where $\aa_j\in S$ is an angle such that $\WWW(\aa_j)= t_j$.
The second identity holds again by rotation invariance.
Replacing in expression (\ref{eert}), 
\begin{align}
\label{27}
\alpha^{-1}=\sum_{i=0}^k \sum_{j=i}^{k} \frac{r_i}{r_k} \frac{K_j}{K_k}
\frac{|\{\sigma\insv : \sigma\cgeq \tilde\om, \sigma_y=0,\sigma_x=\aa_j\}|}
{|\{\sigma\insv : \sigma\cgeq \tilde\om, \sigma_y=0,\sigma_x=0\}|}.
\end{align}
By Lemma \ref{combinatorics} applied to $U=\{y\}$ we get
\begin{align}
  \alpha^{-1}\; \le\; \ds{\sum_{i=0}^k \sum_{j=i}^{k} \frac{r_i}{r_k} \frac{K_j}{K_k}} \;=\; \ds{\sum_{j=0}^k\frac{t_j}{r_k}\frac{K_j}{K_k}},
\end{align}
since $t_j=\sum_{i=0}^j r_i$.
From (\ref{199}), we conclude that the stochastic domination $B_\rho\leq_{st}\phi$ will follow for $\beta$ satisfying
\begin{align}
\rho\;\le\; \varphi(\beta)\;:=\; \left(\sum_{j=0}^k\frac{t_j}{r_k}\frac{K_j}{K_k}\right)^{-1}.
\end{align}
The function $\varphi$ is increasing. Indeed, for each $j$, $\frac{r_k}{t_j}$ is of the form $e^{\beta A}(1-e^{-\beta B})$ with $A$ and $B$ positive numbers, and hence  increasing.
On the other hand $\lim_{\beta\-7\infty}r_k=1$ and $\lim_{\beta\-7\infty}t_i=0$ for $i<k$; as a consequence, $\lim_{\beta\-7\infty}\varphi(\beta)=1$.
Finally, $\lim_{\beta \downarrow 0}r_k=0$ and $\lim_{\beta \downarrow 0}t_i=1$ for every $i$, so $\lim_{\beta \downarrow 0}\varphi(\beta)=0$. See Figure 1 for the graph of $\varphi$ when 
$q=4$.
In particular, $\varphi$ is injective and its inverse $\varphi^{-1}:(0,1)\-7 (0,\infty)$ is well defined.
We conclude that
if $\beta_0=\varphi^{-1}(\rho)$, then equation \eqref{199} holds for $\beta\ge \beta_0$.
\qed


\begin{figure}
\centering
\begin{tikzpicture}[thick,yscale=0.6, xscale=1]
  \draw[<->] (-1,0) -- (7,0)  {};
  \draw[<->] (0,-1) -- (0,7)  {};
  \draw[-] [line width=1.2pt] (-0.09,6) -- (0.09,6)  {};
  \draw[-] [line width=0.3pt] (0.09,2.5) -- (1.69,2.5) [dash pattern=on 1pt off 1pt,domain=-4:4]   ;  
      \draw[-] [line width=0.3pt] (1.69,0) -- (1.69,2.5) [dash pattern=on 1pt off 1pt,domain=-4:4]   ;  
    \draw[-] [line width=1.2pt] (-0.09,2.5) -- (0.09,2.5)  {};
  \draw[-] [line width=0.3pt] (0.12,6) -- (7,6) [dash pattern=on 1pt off 1pt,domain=-4:4]   {};  
  \draw[xscale=1.4,yscale=6,domain=-0:5,smooth,variable=\x] plot ({\x},{ (1-exp(-\x))/(exp(-2*\x)+2*exp(-\x)+1) });
  \draw  (7,0.5) node {$\beta$};
    \draw  (-0.4,2.5) node {$\rho$};
    \draw  (0.56,7) node {$\varphi(\beta)$};
    \draw[-] [line width=1.2pt] (1.69,-0.13) -- (1.69,0.13)  {};
    \draw  (1.7,-0.6) node {$\beta_0(\rho)$};
    \draw  (-0.3,6) node {$1$};
\end{tikzpicture}
\caption{} \label{figure1}
\end{figure}

\paragraph{Proof of Lemma \ref{combinatorics}.}
The case $x\in U$ is trivial, so let us suppose $x\in V\wo U$. If $|U|>1$ the model can be reduced to the case $|U|=1$ by identifying all  vertices in $U$.
We may then suppose $U=\{y\}$ for some $y\neq x$.

Let
\begin{align}
L_\om (a):= \{\sigma\insv :\sigma\cgeq\om, \sigma_y=0, \sigma_x=a\}.
\end{align}
We will construct an injection $\Phi:L_\om(a)\hookrightarrow L_\om(0)$. Here is a brief description of the procedure.
Fix $a\in S$ and consider the reflection $R: S\-7  S$ with respect to the line $\ell$ at angle $a/2$ with the horizontal axis (see Figure \ref{figure2}), that is, $Rb=a-b \ \mbox{mod} \ 2\pi$.
Clearly, $R(a)=0$.
We progressively transform an initial configuration $\sigma\in L_\om(a)$ into
a configuration $\sigma' \in L_\om(0)$ .  The first step is to modify $\sigma$
by applying the reflection $R$ to the spin at the vertex $x$.  The resulting
configuration may present incompatibilities with respect to $\om$ and, if it
does, they will appear at edges $\{\pic{ux}\}_{u\in V}$.  If this is the
case, we modify the configuration by applying the transformation $R$ to the
spins of the conflicting vertices.
 We obtain  a
configuration without incompatibilities in the edges having one endpoint at $x$, but we might have created new incompatibilities at a second level of
edges, that is, edges with one endpoint at a vertex that is a neighbour of $x$. We
solve this by applying $R$ once more to the spins of the new conflicting
vertices, and keep repeating the procedure until there are no more
incompatibilities.
We need to
show that the  resulting configuration $\sigma'$ belongs to $L_\om(0)$, and that the construction is indeed
injective. 
The most delicate part is to prove that this process stops before
reaching the vertex $y$.

\begin{figure}
\centering
\begin{tikzpicture}[line cap=round,line join=round,>=triangle 45,x=1cm,y=1cm]
\path [fill=gray!30!] (-4,-1.5) -- (4,1.5) -- (4,3) -- (-4,3);
\clip(-4,-3) rectangle (4,3);
\draw(0,0) circle (2cm);
\draw [domain=-4:4] plot(\x,{(-0-0*\x)/1});
\draw  (0,-4) -- (0,4);
\draw [line width=1pt,domain=-4:4] plot(\x,{(-0--2.12*\x)/5.66});
\draw [dash pattern=on 1pt off 1pt,domain=-4:4] plot(\x,{(-8.05--5.66*\x)/-2.12});
\draw (0.81,0.67)-- (1.13,0.78);
\draw (0.93,0.35)-- (0.81,0.67);
\draw (3.1,1.75) node[anchor=north west] {$\ell$};
\draw (0.4,-2.12) node[anchor=north west] {$\mbox{Hem}(0)$};
\draw (-2,2.7) node[anchor=north west] {$\mbox{Hem}(a)$};
\draw[fill]   (0,0) circle (1.2pt);
\draw[fill]   (0.72,1.86) circle (1.2pt);
\draw  (0.87,2.05) node {$b$};
\draw[fill]   (1.77,-0.93) circle (1.2pt);
\draw  (2.05,-0.9) node {$Rb$};
\draw[fill]   (1.87,0.7) circle (1.2pt);
\draw  (2.2,0.6) node {\scriptsize $a/2$};
\draw[fill]   (2,0) circle (1.2pt);
\draw  (2.13,-0.12) node {\scriptsize $0$};
\draw[fill]   (-2,0) circle (1.2pt);
\draw  (-2.15,-0.12) node {\scriptsize $\pi$};
\draw[fill]   (1.5034,1.309) circle (1.2pt);
\draw  (1.57,1.48) node {\scriptsize $a$};
\end{tikzpicture}
\caption{} \label{figure2}
\end{figure}

It suffices to prove the result when
$\aa\neq 0$, which we assume from now on.  We may also assume that $\aa\in(0,\pi]$, as the other case is symmetric.  As before, the
boundary $\partial V'$ of a vertex set $V'\subset V$ denotes the set of vertices
$u\in V\wo V'$ such that $\pic{uv}\in \EE$ for some $v\in V'$.

Let now $\sigma\in L_\om(\aa)$.
Define a sequence of sets $A_0\subset A_1\subset\dots\subset V$ associated to $\sigma$ by
$A_0:=\{x\}$
and, for $n\ge 0$, 
\begin{align}
  A_{n+1}:= A_n\cup\llave{u\in \partial A_n: \WWW(\sigma_u-R\sigma_v)<
   \om_{\pic{uv}} \ \mbox{for some} \ v\in A_n}.
\end{align} 
At  each step, $A_{n+1}\setminus A_n$ consists of those vertices where new incompatibilities would arise when applying the reflection to $A_n$. Let
\begin{align}
A:=\bigcup_{n\geq 0} A_n.
\end{align}
Define the function $\Phi:\LLL_\om(\aa)\to S^{V}$ by
\begin{align}
  (\Phi\sigma)_{u} := 
\left\{ 
\begin{array}{ll}
  R\sigma_{u}&\hbox{if } u\in A\\
\sigma_{u}&\hbox{if } u\notin A\,.
\end{array}
\right.
\end{align}
We now show that \emph{i)} the image of $\Phi$ is contained in $L_\om(0)$  and  \emph{ii)} that $\Phi:\LLL_{\om}(\aa)
\rightarrow \LLL_{\om}(0)$ is an injection.

\noindent \emph{i) }$\Phi\sigma\in L_\om(0)$. In order to  prove that $\Phi\sigma\cgeq \om$,
 we need to show that 
 \begin{equation}
 \label{condition}
 \WWW\big((\Phi\sigma)_u-(\Phi\sigma)_v\big)\geq
\om_{\pic{uv}}
\end{equation} 
for any $\pic{uv}\in \EE$.
The cases $\llave{u,v}\subset A$ or $\llave{u,v}\subset
A^c$ are trivial.
If $u\notin A$ and $v\in A$, condition \eqref{condition} reads $\WWW(\sigma_u-R\sigma_v)\geq
\om_{\pic{uv}}$, which must hold;
otherwise $u$ would have belonged to $A$ in the first place.

It remains to prove that $(\Phi\sigma)_y=0$, which follows if we show that $y\notin A$.
The line $\ell$ (see Figure \ref{figure2}) separates the two open hemispheres $\mbox{Hem}(0)$ and $\mbox{Hem}(a)$ defined by
\begin{align*}
&\mbox{Hem}(0):=\llave{b\in S: \sin(b-a/2)<0}
\\
&\mbox{Hem}(a):=\llave{b\in S: \sin(b-a/2)>0}.
\end{align*}
Since $0\in \Hem(0)$, it is enough to prove that $\sigma_u\in \Hem(a)$ for every $u\in A\wo \{x\}$.
We proceed by induction.
If $A_1\neq \emptyset$, let $u\in A_1\wo\{x\}$ with $\sigma_u=b$.
By the definition of $A_1$, we have $\WWW(b-0)<\om_{\pic{ux}}\leq \WWW(b-a)$, where the inequality follows from the fact that $\sigma\cgeq\om$.
Now, $\WWW(b-0)< \WWW(b-a)$ is equivalent to
$\cos(b)< \cos(b-a)$.
But
\begin{align}
\label{first}
\cos(b)< \cos(b-a)\quad
\Longleftrightarrow \ \  \cos\pare{ b -\frac{ a }{2}+\frac{ a }{2}}< \cos\pare{ b -\frac{ a }{2}-\frac{ a }{2}}
\end{align}
\begin{align}
\nonumber
&\Longleftrightarrow \ \ \cos\pare{ b -\frac{ a }{2}}\cos\pare{\frac{ a }{2}}
-\sin\pare{ b -\frac{ a }{2}}\sin\pare{\frac{ a }{2}}
\\[0.5cm]
&\hspace{3.0cm}< \cos\pare{ b -\frac{ a }{2}}\cos\pare{\frac{ a }{2}}
+\sin\pare{ b -\frac{ a }{2}}\sin\pare{\frac{ a }{2}}
\end{align}
\begin{align}
\label{last}
\Longleftrightarrow\quad  0<2\sin\pare{ b -\frac{ a }{2}}\sin\pare{\frac{ a }{2}}\quad
\Longleftrightarrow \quad  0<\sin\pare{ b -\frac{ a }{2}}\quad
\Longleftrightarrow \quad   b \in \Hem(a),
\end{align}
and the claim holds for $A_1$.
Suppose now that $\sigma_u\in\Hem(a)$, that is
\begin{align}
\sin\pare{\sigma_u-\frac{ a }{2}}> 0,
\end{align}
for every $u\in A_n$. If $A_{n+1}\neq \emptyset$, let $v\in A_{n+1}$ and $w\in A_{n}$ be such that  $\WWW(\sigma_v-R\sigma_w)<\WWW(\sigma_v-\sigma_w)$, which is equivalent to $\cos\big(\sigma_v-(a-\sigma_w)\big)<\cos(\sigma_v-\sigma_w)$.
By the inductive hypothesis $\sin\pare{\sigma_w-\frac{a}{2}}>0$.
An argument similar to the one leading from \eqref{first} to \eqref{last} yields
\begin{eqnarray*}
0<2\sin\pare{\sigma_v-\frac{a}{2}}\sin\pare{\sigma_w-\frac{a}{2}},
\end{eqnarray*}
and then $0<\sin\pare{\sigma_v-\frac{a}{2}}$, i.e. $\sigma_v \in \mbox{Hem}(a)$. This completes the induction.

\noindent \emph{ii)} $\Phi$ is injective.  Let $\sigma,\sigma'\in\LLL_{\om}(0,\aa)$ be
two different configurations and denote by $A,A_1,A_2,\ldots$ and $A',A_1',A_2',\ldots$
their associated incompatibility  sets.  If $A=A'$, we are done because $R$ is
injective. Suppose $A\neq A'$ and let
\begin{align*}
n=\min\llave{j\geq 1: A_j\neq A_j'};
\end{align*}
so that in particular $A_{n-1}=A_{n-1}'$.  If there is a vertex $u\in A_{n-1}$ such that
$\sigma_u\neq \sigma'_u$, we are done.  Suppose
$\sigma\stackrel{A_{n-1}}{=}\sigma'$.  Without loss of generality, let us take
$u\in A_n\wo A_n'$.
We claim that $(\Phi \sigma)_u\neq (\Phi \sigma')_u$.
We know that $\sigma_u\neq\sigma'_u$, as otherwise we would have $u\in A_n'$.
If $u\in A'$ we are done.  Suppose then  $u\notin A'$.  Let $v\in A_{n-1}$ be
such that $ \WWW(\sigma_u-R\sigma_v)<\om_{\pic{uv}}$.  Using that $\WWW(\aa'-\bb')=\WWW(R\aa'-R\bb')$ for any
$\aa',\bb'\in S$, that $R^2$ is the identity and that
$(\Phi\sigma)_u=R\sigma_u$, we have $ \WWW(\sigma_u-R\sigma_v)=
\WWW\big((\Phi\sigma)_u-\sigma_v\big)$, and then
\begin{align}\label{pplloo}
 \WWW\big((\Phi\sigma)_u-\sigma_v\big)<\om_{\pic{uv}}.
\end{align}
On the other hand, since $\sigma'\cgeq \om$, we have
$\WWW(\sigma'_u-\sigma'_v)\geq \om_{\pic{uv}}$.
But
$\WWW(\sigma'_u-\sigma'_v)= \WWW\big((\Phi \sigma')_u-\sigma_v\big)$ because $u\notin
A'$ and $\sigma\stackrel{A_{n-1}}{=}\sigma'$, and hence
\begin{align}\label{pplloo1}
 \WWW\big((\Phi \sigma')_u-\sigma_v\big)\geq \om_{\pic{uv}}.
\end{align}
From inequalities (\ref{pplloo}) and (\ref{pplloo1}) we obtain $(\Phi \sigma)_u\neq (\Phi \sigma')_u$, as claimed.
\qed

\section{Appendix}
\label{appx}

\paragraph{Proof of Theorem \ref{edso}: The Edwards-Sokal random-cluster representation.}
The first step is to write the density of $\mu$ with respect to $\hat\mu$:
\begin{align}
\label{jugo}
\mu(\sigma)=\hat\mu( \sigma) \,\frac{1}{Z'} \prod_{\pic{xy}\in
  \EE} \WWW(\sigma_x-\sigma_y),
\end{align}
with $Z'$ the normalizing constant in \eqref{Q}.
Since $\WWW(\sigma_x-\sigma_y)=\theta\big(t\in {\cal I}:
  t\leq \WWW(\sigma_x-\sigma_y)\big)$, the weight of a spin configuration
can be realized as the probability of a related event on the associated edge
set:
\begin{align}
  \prod_{\pic{xy}\in \EE} \WWW(\sigma_x-\sigma_y)
  \;=\;\hat\phi     \Big( \bigcap_{\pic{xy}\in \EE}\big\{\om\in {\cal I}^{\EE}: \om_{\pic{xy}}\leq  \WWW(\sigma_x-\sigma_y)\big\} \Big) \nn 
    \;=\;\hat\phi\big(\om: \om \cleq
    \sigma\big). \label {a1}
\end{align} 
Here is where the definition of compatibility appears naturally.
Inserting \eqref{a1} in \eqref{jugo}, we get
\begin{align}
  \mu(\sigma)
  &=\sum_{\om\in {\cal I}^{\EE}}\frac{1}{
    Z'}  \11 [\om\cleq \sigma] \ \hat\phi(\om) \
  \hat\mu(\sigma)=
  \sum_{\om\in {\cal I}^{\EE}}Q\big((\om,\sigma)\big).
\end{align}
Hence, $\mu$ is the second marginal of $\QQ$.
Adding over all the possible vertex-configurations, it is easy to see that
$\phi$ is its first marginal.
\qed

\paragraph{Asymptotics for $\beta_0$.}

The threshold $\beta_0$ introduced in Theorem \ref{diluted} is $\beta_0=\varphi^{-1}(\rho)$, where $\varphi:(0,\infty)\-7(0,1)$ is the function defined in the proof of Theorem \ref{domination}, and $\rho$ is the parameter defined in the proof of Theorem \ref{diluted}, such that $\rho>\frac{p_c}{p}$. Since $\varphi^{-1}$ is increasing, we can take the infimum over 
$\rho$ to optimize $\beta_0=\varphi^{-1}\big(  \frac{p_c}{p} \big)$.

For any fixed $\beta >0$ we have that  $\lim_{q\-7\infty}\varphi(\beta)=0$.
Indeed, $\lim_{q\-7\infty}r_k=0$ and $t_i\frac{K_i}{K_k}$ is bounded away from zero uniformly in $q$.
As a consequence, for every fixed $\tilde{p}\in (0,1)$, $\lim_{q\-7\infty}\varphi^{-1}(\tilde{p})=\infty$.
We conclude that our method is not informative as a discretization of the  $XY$ model, that is, when the  number of spins goes to infinity.

Note that
\begin{align}
\varphi(\beta)^{-1}=\frac{1}{r_k}+\sum_{i=0}^{k-1}\frac{t_i}{r_k}\frac{K_i}{K_k}
\leq \frac{1}{r_k}+\sum_{i=0}^{k-1}\frac{t_{k-1}}{r_k}\ 2
\leq \frac{1}{r_k}+\frac{t_{k-1}}{r_k}\ q,
 \end{align} 
so that
\begin{align}
\varphi(\beta)\geq \frac{r_k}{1+qt_{k-1}}.
\end{align}
Then $\beta_0$ is bounded above by the solution to the  equation
\begin{align}
\frac{p_c}{p}=\frac{r_k}{1+qt_{k-1}}.
\end{align}
Using that $r_k=1-t_{k-1}$ and that $t_{k-1}=e^{-\beta \paree{1-\cos\paree{\frac{2\pi}{q}}}}$, this solution can be explicitly computed as
\begin{align}
\frac{\log\paree{\frac{p+qp_c}{p-p_c}}}{1-\cos\paree{\frac{2\pi}{q}}}.
\end{align}
If we fix  $p$ and $d$, this expression is
 of order $q^2\log(q)$ as $q\-7 \infty$, the same order given by  Pirogov-Sinai theory and reflection positivity in the $2$-dimensional homogeneous case.
If we fix $p$ and $q$, it is of order 
\begin{align}
\frac{\log\paree{1+\frac{1}{d}}}{1-\cos\paree{\frac{2\pi}{q}}}
\end{align}
as $d\-7 \infty$, taking into account that $p_c\sim \frac{1}{2d}$.
In particular, $\beta_0\-7 0$ as $d\-7 \infty$.

\section*{Acknowledgments}
We are grateful to both referees for several comments that helped us give a clear focus to the article.  
Nahuel Soprano-Loto would like to thank Roberto Fern\'andez for many useful discussions and his warm welcome to the University of Utrecht, a visit funded by grant 2010-MINCYT-NIO Interacting stochastic systems: fluctuations, hydrodynamics, scaling limits.  This research has been supported by grant PICT 2012-2744 ÒStochastic Processes and Statistical MechanicsÓ.

\bibliographystyle{alpha}
\bibliography{dilute-clock}

\begin{thebibliography}{ACCN87}

\bibitem[ACCN87]{ACCN87}
M.~Aizenman, J.~T. Chayes, L.~Chayes, and C.~M. Newman.
\newblock The phase boundary in dilute and random {I}sing and {P}otts
  ferromagnets.
\newblock {\em J. Phys. A}, 20(5):L313--L318, 1987.

\bibitem[Bis09]{biskup}
Marek Biskup.
\newblock Reflection positivity and phase transitions in lattice spin models.
\newblock In {\em Methods of contemporary mathematical statistical physics},
  volume 1970 of {\em Lecture Notes in Math.}, pages 1--86. Springer, Berlin,
  2009.

\bibitem[Dob68]{dobrushin}
R.~L. Dobrushin.
\newblock The description of a random field by means of conditional
  probabilities and conditions of its regularity.
\newblock {\em Theor. Prob. Appl.}, 13:197--224, 1968.

\bibitem[Dub11]{D11}
Julien Dub{\'e}dat.
\newblock Topics on abelian spin models and related problems.
\newblock {\em Probab. Surv.}, 8:374--402, 2011.

\bibitem[ES88]{es}
Robert~G. Edwards and Alan~D. Sokal.
\newblock Generalization of the {F}ortuin-{K}asteleyn-{S}wendsen-{W}ang
  representation and {M}onte {C}arlo algorithm.
\newblock {\em Phys. Rev. D (3)}, 38(6):2009--2012, 1988.

\bibitem[FILS78]{FILS78}
J{\"u}rg Fr{\"o}hlich, Robert Israel, Elliott~H. Lieb, and Barry Simon.
\newblock Phase transitions and reflection positivity. {I}. {G}eneral theory
  and long range lattice models.
\newblock {\em Comm. Math. Phys.}, 62(1):1--34, 1978.

\bibitem[FSL]{F-SL}
Roberto Fern\'andez and Nahuel Soprano-Loto.
\newblock {\em In preparation}.

\bibitem[GHM01]{GHM01}
Hans-Otto Georgii, Olle H{\"a}ggstr{\"o}m, and Christian Maes.
\newblock The random geometry of equilibrium phases.
\newblock In {\em Phase transitions and critical phenomena, {V}ol. 18},
  volume~18 of {\em Phase Transit. Crit. Phenom.}, pages 1--142. Academic
  Press, San Diego, CA, 2001.

\bibitem[Gri06]{G06}
Geoffrey Grimmett.
\newblock {\em The random-cluster model}, volume 333 of {\em Grundlehren der
  Mathematischen Wissenschaften [Fundamental Principles of Mathematical
  Sciences]}.
\newblock Springer-Verlag, Berlin, 2006.

\bibitem[KF69]{KF69}
P.~W. Kasteleyn and C.~M. Fortuin.
\newblock Phase transitions in lattice systems with random local properties.
\newblock {\em Journal of the Physical Society of Japan}, 26:11--14, 1969.

\bibitem[MS11]{MS11}
Christian Maes and Senya Shlosman.
\newblock Rotating states in driven clock- and {XY}-models.
\newblock {\em J. Stat. Phys.}, 144(6):1238--1246, 2011.

\bibitem[PS75]{PS75}
S.~A. Pirogov and Ja.~G. Sina{\u\i}.
\newblock Phase diagrams of classical lattice systems.
\newblock {\em Teoret. Mat. Fiz.}, 25(3):358--369, 1975.

\bibitem[vdBM94]{vdbm}
J.~van~den Berg and C.~Maes.
\newblock Disagreement percolation in the study of {M}arkov fields.
\newblock {\em Ann. Probab.}, 22(2):749--763, 1994.

\bibitem[vEKO11]{vEKO11}
Aernout C.~D. van Enter, Christof K{\"u}lske, and Alex~A. Opoku.
\newblock Discrete approximations to vector spin models.
\newblock {\em J. Phys. A}, 44(47):475002, 11, 2011.

\end{thebibliography}

\end{document}